\documentclass[a4paper,12pt,reqno]{amsart}

\usepackage{amsfonts,amssymb,palatino,euler}
\usepackage{stmaryrd}
\usepackage[left=3cm,right=3cm,top=3cm,bottom=3cm]{geometry}
\usepackage{cite}%
\usepackage[colorlinks=true]{hyperref}

\title{Overgroups of elementary groups in polyvector representations}

\author{Roman~Lubkov}
\thanks{The research was supported by the RFBR grant 19-31-90072 and by the Leonhard Euler International Mathematical Institute, agreement no. 075–15–2019–1620}
\address{St.~Petersburg Department of V.\,A.\,Steklov Institute of Mathematics of the Russian Academy of Sciences}
\email{RomanLubkov@yandex.ru}



\let\opn\operatorname
\DeclareMathOperator{\E}{E}
\DeclareMathOperator{\GL}{GL}
\DeclareMathOperator{\sign}{sgn}
\DeclareMathOperator{\Tran}{Tran}
\DeclareMathOperator{\Lie}{Lie}

\DeclareMathOperator{\EO}{EO}
\DeclareMathOperator{\GO}{GO}
\DeclareMathOperator{\GSp}{GSp}
\DeclareMathOperator{\SL}{SL}
\DeclareMathOperator{\level}{lev}

\renewcommand{\trianglelefteq}{\trianglelefteqslant}
\renewcommand{\leq}{\leqslant}
\renewcommand{\geq}{\geqslant}

\newcommand{\bw}[1]{\mathord{\raisebox{2pt}
{\hbox{$\scriptstyle{\bigwedge^{\!#1}}$}}}}
\newcommand\blank{\mathord{\hbox to 1.5ex{\hrulefill}}\,}

\theoremstyle{plain}
\newtheorem{theorem}{Theorem}
\newtheorem{problem}{Problem}
\newtheorem{prop}[theorem]{Proposition}
\newtheorem{lemma}[theorem]{Lemma}
\newtheorem{corollary}[theorem]{Corollary}
\theoremstyle{remark}
\newtheorem*{remark}{Remark}

\newenvironment{thmbis}
  {\addtocounter{theorem}{-1}%
   \begin{theorem}}
  {\end{theorem}}

\begin{document}

\begin{abstract}
We initiate the study of subgroups $H$ of the general linear group $\GL_{\binom{n}{m}}(R)$ over a commutative ring $R$ that contain the $m$-th exterior power of an elementary group $\bw{m}\E_n(R)$. Each such group $H$ corresponds to a uniquely defined level $(A_0,\dots,A_{m-1})$, where $A_0,\dots,A_{m-1}$ are ideals of $R$ with certain relations. In the crucial case of the exterior squares, we state the subgroup lattice to be standard. In other words, for $\bw{2}\E_n(R)$ all intermediate subgroups $H$ are parametrized by a single ideal of the ring $R$. Moreover, we characterize $\bw{m}\GL_n(R)$ as the stabilizer of a system of invariant forms. This result is classically known for algebraically closed fields, here we prove the corresponding group scheme to be smooth over $\mathbb{Z}$. So the last result holds over arbitrary commutative rings.
\end{abstract}

\maketitle

\section*{Introduction}\label{sec:intro}

Over the last 40 years one of the central problems in the theory of finite groups was description of their maximal subgroups. In 1984, Michael Aschbacher established the Subgroup Structure Theorem~\cite{AschbacherClasses}. It defines eight explicitly described classes $\mathcal{C}_1$--$\mathcal{C}_8$ and an exceptional class $\mathcal{S}$ for all maximal subgroups of a finite classical group. Later this theorem was reinterpreted several times. For instance, in 1998 Martin Liebeck and Gary Seitz revised this theorem using the theory of algebraic groups\cite{LiebSeitStrucClass}.

For some special cases of fields subgroups of groups from Aschbacher classes have been studied intensely. Over a finite field maximality of subgroups from Aschbacher classes was obtained by Peter Kleidman and Martin Liebeck in the book~\cite{KleidLiebSubStruct}. The same problem for arbitrary fields was solved in a series of publications by Oliver King, Roger Dye, Shang Zhi Li, and others. Occasionally, these subgroups are not maximal, and then the issue was to describe their overgroups. We recommend the surveys~\cite{VavSbgs,VavStepSurvey}, which contain necessary preliminaries, the complete history, and many further related references.

Until recently, over arbitrary rings little was known about the description of overgroups. The first attempts to transfer such results were initiated by Zenon Borevich and Nikolai Vavilov for the Aschbacher classes $\mathcal{C}_1+\mathcal{C}_2$. Later for other classes overgroups were described by the St.~Petersburg algebraic school, Shang Zhi Li and others.

In the present paper, we consider the elementary group in the $m$-th vector representation $\bw{m}\E_n(R)$ over a commutative ring $R$. This case pertains to the exceptional Aschbacher class $\mathcal{S}$. There were extremely few results for this class over arbitrary rings. 

Let us mention the known results on description of overgroups for $\bw{m}\E_n(R)$. Over finite fields Bruce Cooperstein proved maximality of the normalizer $N\bigl(\bw{2}\E_n(K)\bigr)$ in the general linear group $\GL_{\binom{n}{2}}(K)$~\cite{CoopersteinStructure}. For algebraically closed fields description of overgroups of $\bw{m}\E_n(K)$ follows from the results of Gary Seitz on maximal subgroups of classical algebraic groups, for instance, see~\cite{SeitzMaxSub,burness_testerman_CorrSeitz}.

For a commutative ring $R$ the elementary group $\E_n(R)$ is generated by the elementary transvections $t_{i,j}(\xi)$, where $1\leq i\neq j\leq n, \xi\in R$. Thus $\bw{m}\E_n(R)$ is a subgroup of $\GL_{\binom{n}{m}}(R)$ generated by the images of all elementary transvections under the representation $\bw{m}$, see \S\ref{sec:definition} for the precise definitions. In the sequel, we assume that $n\geq 4$ and $m<n$, otherwise the problem degenerates: $\bw{m}\E_3(R)\cong\E_3(R)$ and $\bw{n}\E_n(R)\cong\E_n(R)$.

We initiate the classification of overgroups $H$ such that
$$\bw{m}\E_n(R)\leq H \leq \GL_{\binom{n}{m}}(R).$$

The conjectural answer is the \textit{standard description of overgroups}, which can be formulated as follows. Each intermediate subgroup $H$ corresponds to a uniquely defined level of this subgroup $\level(H)$. In general, the level is an $m$-tuple of ideals $(A_0,\dots,A_{m-1})$ in the ring $R$ with certain relations. 

The plausibility of the results for $\bw{m}\E_n(R)$ seems obvious due to the close connection of the $m$-th vector representation and the tensor product of elementary linear groups. Recall that all overgroups of $\E_k(R)\otimes\E_l(R)$ are parametrized by a triple of ideals $(A,B,C)$ in the ring $R$, see~\cite{AnaVavSinI}. Moreover, in the important particular case $m=2$, $n=4$ the problem was completely solved in the paper~\cite{VP-EOeven} devoted to elementary orthogonal groups\footnote{$\bw{2}\E_4(R)\cong\EO_6(R)$}. For $\bw{2}\E_4(R)$ the level of an overgroup consists of one ideal $A$ under the assumption that $2$ is invertible in $R$. And, unlike the symplectic case, this restriction cannot be reduced.

However for an arbitrary $m$-th vector representation one ideal is not enough to describe all overgroups of $\bw{m}\E_n(R)$ with any such simplifying assumptions as $R$ being Noetherian or $2$ being invertible in $R$. Therefore initially one has to analyze the case of the bivector representation, for which one ideal is [almost] sufficient.

In the present paper, we outline the general scheme of classification of overgroups for $\bw{m}\E_n(R)$. Classically, it consists of the following steps.

$\bullet$ To calculate the level and the normalizer of connected (i.\,e., perfect) intermediate subgroups $H$. In \S\ref{sec:LevelComput} we associate each such a subgroup with the \textit{lower level} of $H$.

$\bullet$ To extract an elementary transvection from $H$. Unlike the classical groups, for $\bw{m}E_n(R)$ the standard methods of extraction are unaviable. Thus we use an alternative technique developed by Alexey Stepanov and the author. For this purpose in \S\ref{sec:InvForms} we construct a system of invariant forms for $\bw{m}\GL_n(\blank)$.

$\bullet$ And the last step is to prove the inclusion in the normalizer. Using the method from the previous item, in \S\ref{sec:extraction} we state that the subgroup lattice is standard. In other words, we establish the required inclusion.

\section{Exterior powers of elementary groups}\label{sec:definition}

In this section, we define an exterior power of the elementary group in detail. By $[n]$ we denote the set $\{1,2,\dots, n\}$ and by $\bw{m}[n]$ we denote an exterior power of the set $[n]$. Elements of $\bw{m}[n]$ are ordered subsets $I\subseteq [n]$ of cardinality $m$ without repeating entries:
$$\bw{m}[n] = \{(i_{1},\dots,i_{m})\mid 1\leq i_1<i_2<\dots<i_m\leq n \}.$$

Let $R$ be a commutative ring and let $R^n$ be the right free $R$-module with the standard basis $\{e_1,\dots,e_n\}$, $n\geq 3$. \textit{By $N$ we denote the binomial coefficient $\binom{n}{m}$}. $\bw{m}R^n$ is the free module of rank $N=\binom{n}{m}$ with the basis $e_{i_1}\wedge\dots\wedge e_{i_m}$, where $1\leq i_1<\dots<i_m\leq n$. Products $e_{i_1}\wedge\dots\wedge e_{i_m}$ are defined for any set $i_1,\dots,i_m$ as $e_{\sigma(i_{1})}\wedge\dots\wedge e_{\sigma(i_{m})} = \sign(\sigma)\, e_{i_1} \wedge \dots \wedge e_{i_{m}}$ for any permutation $\sigma$ in the permutation group $S_{m}$.

For every $m\leq n$ define $\bw{m}$ as a homomorphism from $\GL_n(R)$ into $\GL_N(R)$ by
$$\bw{m}(g)(e_{i_1}\wedge\dots\wedge e_{i_m}):=(ge_{i_1})\wedge\dots\wedge (ge_{i_m})$$
for every $e_{i_1},\dots,e_{i_m}\in R^n$. Thus $\bw{m}$ is a representation of the group $\GL_n(R)$. It is called the $m$-th \textit{vector representation} or the $m$-th \textit{fundamental representation}. $\bw{m}\GL_n(R)$ is called the $m$-th exterior power of the general linear group. The \textit{$[$absolute$]$ elementary group} $\E_n(R)$ is a subgroup of $\GL_n(R)$ generated by all elementary transvections $t_{i,j}(\xi)=e+\xi e_{i,j}$, where $1\leq i\neq j\leq n$, $\xi\in R$. Therefore the exterior power of the elementary group is well defined.

For arbitrary rings there is a difference between $\bw{m} \bigl(\GL_{n}(R)\bigr)$ and $\bw{m}\GL_{n}(R)$. The first group is a set-theoretic image of the [abstract] group $\GL_{n}(R)$ under the Cauchy--Binet homomorphism $\bw{m}\colon\GL_{n}(R)\longrightarrow \GL_N(R)$, whereas $\bw{m}\GL_{n}(\blank)$ is the \textbf{categorical} image a group scheme $\GL_{n}(\blank)$. For rings $\bw{m}\GL_n(R)$ is strictly larger than $\bw{m}\bigl(\GL_{n}(R)\bigr)$, see~\cite{VavPere} for detail.

\section{Level computation and level reduction}\label{sec:LevelComput}

Let $H$ be an overgroup of the exterior power of the elementary group $\bw{m}\E_n(R)$:
$$\bw{m}\E_n(R)\leq H \leq \GL_N(R).$$
And for any indices $I, J \in \bw{m}[n]$ let $A_{I,J}$ be subsets of the ring $R$:
$$A_{I,J}:=\{\xi \in R \;|\; t_{I,J}(\xi) \in H \}$$
As usual, diagonal sets $A_{I,I}$ equal the whole ring $R$ for any index $I\in\bw{m}[n]$. Thus we will construct a $D$-net of ideals for the ring $R$ in terms of Zenon Borevich~\cite{BVnets}. Define a \textit{distance} between indices $I$ and $J$ as the cardinality of the set $I\cap J$:
$$d(I,J)=|I\cap J|.$$
This combinatorial characteristic plays the same role as the distance function $d(\lambda,\mu)$ for roots $\lambda$ and $\mu$ on the weight diagram of a root system. 

It turns out that $A_{I,J}$ are ideals in the ring $R$. Moreover, $A_{I,J}$ \textit{depends not on} $I,J$ themselves, but only on the distance of $(I,J)$.
\begin{prop}\label{prop:Relations}
The ideals $\{A_0,\dots,A_{m-1}\}$ are interrelated as follows:
\begin{gather*}
A_k\leq A_{k+1}, \text{ for } n\geq 3m-2k;\\
A_{0} \geq A_{1} \geq A_{2} \geq \dots \geq A_{m-2} \geq A_{m-1};\\
{\scriptstyle\binom{n-2}{m-1}}\cdot A_{m-2}\leq A_{m-1}.
\end{gather*}
\end{prop}

Note that if $n\geq 3m$, then all ideals coincide. And then the set $A = A_{I,J}$ is called a \textit{level} of an overgroup $H$. Conversely, for $n< 3m$ a level consists of up to $m$ ideals $(A_0,\dots,A_{m-1})$.

To formulate the following theorem define a \textit{$[$relative$]$ elementary group of level $A$}, where $A$ is an ideal in $R$. This group is a normal closure of $\E_n(A)$ in $\E_n(R)$:
$$\E_n(R,A):=\langle t_{i,j}(\xi), 1\leq i\neq j\leq n, \xi\in A\rangle^{\E_n(R)}.$$

\begin{theorem}[Level computation]\label{thm:LevelForm}
Let $R$ be a commutative ring and let $n\geq 3m$. For an arbitrary overgroup $H$ of $\bw{m}\E_n(R)$ there exists a unique maximal ideal $A$ of the ring $R$ such that 
$$\bw{m}\E_n(R)\cdot\E_N(R,A)\leq H.$$
Namely, if a transvection $t_{I,J}(\xi)$ belongs to the group $H$, then $\xi\in A$.
\end{theorem}

In the general case, the level computation is formulated with a set of ideals. An $m$-tuple of ideals
$\mathbb{A}=(A_0,\dots,A_{m-1})$ of the ring $R$ is called \textit{admissible} if $\mathbb{A}$ satisfies the relations in Proposition~\ref{prop:Relations}. Then every admissible $m$-tuple $\mathbb{A}$ corresponds to the group $\E\bw{m}\E_n(R,\mathbb{A}) := \bw{m}\E_n(R)\cdot\E_N(R,\mathbb{A})$. This group is defined as a subgroup generated by $\bw{m}\E_n(R)$ and by all elementary transvections $t_{I,J}(\xi)$, where $\xi\in A_{I,J}$:
$$\E\bw{m}\E_n(R,\mathbb{A})=\bw{m}\E_n(R)\cdot\langle t_{I,J}(\xi),\;\xi\in A_{I,J}\rangle.$$
\begin{thmbis}[Level computation]\label{thm:LevelFormIdeals}
Let $R$ be a commutative ring and let $n\geq 4$. For an arbitrary overgroup $H$ of the group $\bw{m}\E_n(R)$ there exists a net of ideals $\mathbb{A}$ of the ring $R$ such that
$$\bw{m}\E_n(R)\cdot\E_N(R,\mathbb{A})\leq H.$$
Namely, if a transvection $t_{I,J}(\xi)$ belongs to the group $H$, then $\xi\in A_{I,J}$.
\end{thmbis}

Now let $\rho_{A}\colon\GL_{n}(R)\longrightarrow \GL_{n}(R/A)$ be the reduction homomorphism. In the following theorem we describe the normalizer of the lower bound for $H$.
\begin{theorem}[Level reduction]\label{thm:LevelReduction}
Let $n\geq 3m$.
For any ideal $A\trianglelefteq R$, we have
$$N_{\GL_{N}(R)}\bigl(\E\bw{m}\E_n(R,A)\bigr)=\rho_{A}^{-1}\bigl(\bw{m}\GL_{n}(R/A)\bigr).$$
\end{theorem}
Observe that for an $m$-tuple of ideals $\mathbb{A}$ in $R$ it is not obvious how to define congruence subgroups $\rho_{A}^{-1}\bigl(\bw{m}\GL_{n}(R/\mathbb{A})\bigr)$, not to mention the level reduction. The author work in this direction.

\section{Construction of the invariant forms}\label{sec:InvForms}

Our immediate goal is to define $\bw{m}\GL_n(R)$ as a stabilizer of certain invariant forms. Over algebraically closed fields these results are well known. In this section, we construct invariant forms over arbitrary rings.

\subsection{Exterior powers as a stabilizer of invariant forms I}\label{subsec:StabI}

We assume that $2\in R^{*}$ and $n \geq 2m$ due to the isomorphism $\bw{m}V^{*} \cong (\bw{\dim(V)-m}V)^{*}$ for an arbitrary free $R$-module $V$. The following theorem is classically known and can be found in \cite[Chapter 2, Sections 5--7]{DieCarInvTheor}.
\begin{prop}\label{prop:FormsGood}
Let $K$ be an algebraically closed field. Then the group $\bw{m}\GL_{n}(K)$ has an invariant form only in the case $\frac{n}{m} \in \mathbb{N}$, and then the form is unique and equals 
\begin{itemize}
\item $q^m_{[n]}(x) = \sum \sign(I_{1}, \dots, I_{\frac{n}{m}})\; x_{I_{1}}\dots x_{I_{\frac{n}{m}}}$ for even $m$;
\item $q^m_{[n]}(x) = \sum \sign(I_{1}, \dots, I_{\frac{n}{m}})\; x_{I_{1}}\wedge \dots \wedge x_{I_{\frac{n}{m}}}$ for odd $m$,
\end{itemize}
where the sums in the both cases range over all unordered partitions of the set $[n]$ into $m$-element subsets $I_{1}, \dots, I_{\frac{n}{m}}$.
\end{prop}

Thus $\bw{m}\GL_{n}(K)$ is isomorphic to the group of matrices $g\in\GL_N(K)$ for which there is a multiplier $\lambda\in R^*$ such that $q^m_{[n]}(gx)=\lambda(g)q^m_{[n]}(x)$ for all $x\in K^N$, where $K$ is an algebraically closed field. Obviously, $\lambda(g)$ is a one-dimensional representation of the group $\GL_n(K)$. Therefore $\lambda=\det^{\otimes l}\colon g \mapsto \det^{l}(g)$. To evaluate the power of the determinant, we calculate $\lambda(g)$ for the diagonal matrix $d_i(\xi)\in\GL_n(K)$. Thus
$q^m_{[n]}(\bw{m}d_i(\xi)\cdot x)=\xi q^m_{[n]}(x)$. This implies that $\lambda(g)=\det(g)$.

In the sequel, we use the uniform notation $q(x)$ for these forms. This cannot lead to a confusion as we always can distinguish between two different meanings of $q(x)$ by the power $m$. First note that these forms are the only possible ones for the group $\bw{m}\GL_{n}(R)$. Further, the coefficients equals $\pm 1$, so they are defined over $\mathbb{Z}$. Then by direct calculations we get
$$q(\bw{m}g \cdot x) = \det(g) \cdot q(x) \text{ for any }g \in \GL_{n}(R).$$
Thus we can assume these forms to be invariant under the action $\bw{m}\GL_n(R)$, where $R$ is a commutative ring. Observe that it is easy to prove $\bw{m}\E_n(R)$ to preserve these forms.
\begin{prop}\label{prop:FormsInvariantUnderLE}
Let $R$ be an arbitrary commutative ring. If $\frac{n}{m} \in \mathbb{N}$, then the form $q(x)$ is invariant under the action of $\bw{m}\E_n(R)$.
\end{prop}

Suppose that $m$ is even for the sake of brevity. Now for the form $q(x)$ we introduce a \textit{$[$full$]$ polarization}. Put
$k:=\frac{n}{m}\in\mathbb{N}$ and define a $k$-linear form:
$$f^m_{[n]}(x^1,\dots,x^k)=\sum \sign(I_{1}, \dots, I_k)\; x^1_{I_{1}}\dots x^{k}_{I_k},$$
where the sum ranges over all \textit{ordered} partitions of the set $[n]$ into $m$-element subsets. For odd $m$ the form $f^m_{[n]}(x^1,\dots,x^k)$ is set similarly. Again, we denote these forms by the uniform symbol $f(x^1,\dots,x^k)$. We focus on ordered partitions in this sum, unlike unordered ones in $q(x)$.

Let us define a group $G_f(R)$ as the group of linear transformations preserving the form $f$:
$$G_f(R):=\{g\in\GL_N(R)\mid f(gx^1,\dots,gx^k)=f(x^1,\dots,x^k)\text{ for all } x^1,\dots,x^k\in R^N\}.$$
It is an analogue of the Chevalley group for the exterior powers. The extended Chevalley group satisfies similarities of $f$:
\begin{align*}
\overline{G}_f(R)&:=\{g\in\GL_N(R)\mid\text{ there exists } \lambda\in R^* \text{ such that }
\\
&\quad f(gx^1,\dots,gx^k)=\lambda(g)f(x^1,\dots,x^k)\text{ for all } x^1,\dots,x^k\in R^N\}.
\end{align*}
Obviously, the functors $R\mapsto\overline{G}_f(R)$ and $R\mapsto G_f(R)$ define affine group schemes over $\mathbb{Z}$. Thus in the case $k:=\frac{n}{m} \in \mathbb{N}$ we can expect the group $\bw{m}\GL_{n}(R)$ to coincide with the stabilizer of the corresponding form. This is almost true and the following theorem gives a precise answer.
\begin{theorem}\label{thm:StabAndGenerealPower}
Suppose $\frac{n}{m} \in \mathbb{N}$; then $\bw{m}\GL_{n}(R)$ coincides with $\overline{G}_f(R)$ except the case of half dimension. Namely if $n = 2m$, then $\overline{G}_f(R) = \GO_{N}(R)$ or $\GSp_{N}(R)$ depending on the parity of $m$. So in this case $\bw{m}\GL_{n}(R)$ is a subgroup of the orthogonal or the symplectic group respectively;
\end{theorem}
\begin{remark}
If $(n,m) = (4,2)$, then the stabilizer equals $\GO_{6}(R)$. And it also coincides with $\bw{2}\GL_{4}(R)$.
\end{remark}

The proof uses\footnote{Similarly, the statement could be proved using SGA, Exp. VI$\_$b, Cor. 2.6} the Waterhouse lemma, see~\cite[Theorem~$1.6.1$]{WaterhouseDet}. This result reduces the verification of an isomorphism of affine group schemes to the isomorphism of their groups of points over algebraically closed fields and the dual numbers over such fields. In the proof there is two key steps. Step 1: we must check that the theorem holds for any algebraically closed field. And Step 2: the schemes $\overline{G}_f$, $G_f$ must be smooth. This is essentially the same, to evaluate the dimension of the Lie algebras.
\begin{prop}\label{prop:dimLie}
In the non-exceptional case $n\neq 2m$ for any field $K$ the dimension of the Lie algebra $\Lie(\overline{G}_f(K))$ does not exceed $n^2$, whereas the dimension of the Lie algebra $\Lie(G_f(K))$ does not exceed $n^2-1$.
\end{prop}

Consequently, we establish the required isomorphisms.
\begin{theorem}\label{thm:exGLandForms}
If $n\neq 2m$, then there are isomorphisms $\overline{G}_f\cong\bw{m}\GL_n$, $G_f\cong\bw{m}\SL_n$ of affine groups schemes over $\mathbb{Z}$.
\end{theorem}

This theorem guarantees that for arbitrary rings the class of transvections from $\bw{m}\GL_n(R)$ is strictly larger than the images $\bw{m}g$, where $g\in\GL_n(R)$, for detail see~\cite[\S9]{VavPere} for an arbitrary power or~\cite{WaterhouseDet} for the exterior squares. Below we give a similar analysis for the special linear group. The exact sequence of affine group schemes
$$1\longrightarrow\mu_d\longrightarrow\SL_n\longrightarrow\SL_n/\mu_d\longrightarrow 1$$
gives the exact sequence of Galois cohomology
\begin{multline*}
    1\longrightarrow\mu_d(R)\longrightarrow\SL_n(R)\longrightarrow\SL_n/\mu_d(R)\longrightarrow\\
    H^1(R,\mu_d)\longrightarrow H^1(R,\SL_n)\longrightarrow H^1(R,\SL_n/\mu_d),
\end{multline*}
where $d=\gcd(n,m)$. The values of all these cohomology sets are well known, see, for instance~\cite[Chapter~III, \S2]{Knus}.

The set $H^1(R,\mu_d)$ classifies projective modules $P$ of rank $1$ together with the isomorphism $P^{\otimes d}=R$, whereas the map $H^1(R,\mu_d)\longrightarrow H^1(R,\SL_n)$ sends a projective module $P$ to the direct sum $\bigoplus_1^n P=P\oplus\dots\oplus P$. Thus its kernel, as well as the quotient group of $\bw{m}\SL_n(R)$ modulo $\bw{m}\bigl(\SL_n(R)\bigr)$, contains projective modules $P$ of rank $1$ such that $P^{\otimes d}=R$ and $P\oplus\dots\oplus P=R^n$.

\subsection{Exterior powers as a stabilizer of invariant forms II}\label{subsec:StabII}

In the previous subsection, we completely analyzed the case of one invariant form. But if $\frac{n}{m} \not\in\mathbb{N}$, then it turns out that the group $\bw{m}\GL_n(R)$ has an ideal of invariant forms. Let us extend the definition of $q(x)$ from~\S\ref{subsec:StabI}. By default, the form $q(x)$ is considered for the set $[n]=\{1,\dots,n\}$. But in the sequel, we use forms for certain subsets of $[n]$. Let $V\subseteq [n]$ is a $n_1$-subset of $[n]$, where $\frac{n_1}{m}\in\mathbb{N}$. Define a form $q^m_V(x)$ as follows.
\begin{itemize}
\item $q^m_V(x)=\sum \sign(I_{1}, \dots, I_{\frac{n_1}{m}})\; x_{I_{1}}\dots x_{I_{\frac{n_1}{m}}}$ for even $m$;
\item $q^m_V(x)=\sum \sign(I_{1}, \dots, I_{\frac{n_1}{m}})\; x_{I_{1}}\wedge \dots \wedge x_{I_{\frac{n_1}{m}}}$ for odd $m$,
\end{itemize}
where the sums in the both cases range over all unordered partitions of the set $V$ into $m$-element subsets $I_{1}, \dots, I_{\frac{n_1}{m}}$.

As above, $f^m_{V}(x^1,\dots,x^k)$ is the \textit{$[$full$]$ polarization} of $q^m_V(x)$, where $k:=\frac{n_1}{m}$. Furthermore, we ignore the power $m$ in the notation $f^m_{V}(x^1,\dots,x^k)$ and $q^m_V(x)$.

Let divide $n$ by $m$ with the remainder: $n = lm +r$, where $l,r \in \mathbb{N}$. Consider the ideal $Y=Y_{n,m}$ of the ring $\mathbb{Z}[x_I]$ generated by the forms $f_{V}(x^1,\dots,x^k)$ for all possible $m\cdot l$-element subsets $V\subset [n]$.

Suppose that $Y$ is generated by $f_{V_1},\dots,f_{V_p}$, where $p=\binom{n}{ml}$. Then define the extended Chevalley group $\overline{G}_Y(R)$ as the group of linear transformations preserving the ideal $Y$:
\begin{align*}
\overline{G}_Y(R)&:=\{g\in\GL_N(R)\mid\text{ there exist } \lambda_{V_1},\dots,\lambda_{V_p}\in R^*, c(V_k,V_l)\in R \text{ such that }
\\
&\quad f_{V_j}(gx^1,\dots,gx^k)=\lambda_{V_j}(g)f_{V_j}(x^1,\dots,x^k)+\sum\limits_{l\neq j}c({V_j},{V_l})\cdot f_{V_l}(x^1,\dots,x^k)
\\
&\qquad\qquad\qquad\qquad\qquad\qquad\qquad\qquad\text{ for all } 1\leq j\leq p\text{ and } x^1,\dots,x^k\in R^N\}.
\end{align*}

In the rest of the section, we verify that $\bw{m}\GL_n(R)$ is the extended Chevalley group $\overline{G}_Y(R)$ for $m\nmid n$. Again, we use the Waterhouse lemma for this purpose. But now it is not obvious that $\overline{G}_Y$ is a group scheme.
\begin{lemma}
Let $n=ml+r$, where $m,l\in\mathbb{N}$. Then the functor $R\mapsto\overline{G}_Y(R)$ is an affine group scheme over $\mathbb{Z}$.
\end{lemma}

Arguing as in the previous subsection, we must only check that $\overline{G}_Y$ is smooth and the statement holds for any algebraically closed field.
\begin{prop}\label{prop:dimLieIdeal}
For any field $K$ the dimension of the Lie algebra $\Lie(\overline{G}_Y(K))$ does not
exceed $n^2$.
\end{prop}

Consequently, we conclude that $\bw{m}\GL_n(\blank)$ equals the stabilizer of $Y$.
\begin{theorem}\label{thm:exGLandFormsIdeal}
If $n=ml+r$, where $m,l\in\mathbb{N}$, then there is an isomorphism $\overline{G}_Y\cong\bw{m}\GL_n$ of affine group scheme over $\mathbb{Z}$.
\end{theorem}

\section{Subgroup lattice is standard}\label{sec:extraction}

The extraction of a unipotent (an elementary transvection) from an intermediate subgroup $H$ is the key point for the standard description of overgroups. In the previous papers on this problem, e.\,g., joint works of Nikolai Vavilov and Victor Petrov on overgroups of classical groups, the authors extract a nontrivial unipotent from an arbitrary element $a\in G(\Phi,R)\smallsetminus N\bigl(E(\Phi,R)\bigr)$. In fact, for this purpose they extract a number of unipotents and prove that all these unipotents do not vanish under the canonical homomorphism $g\mapsto a$, where $g$ is the generic element of $G(\Phi,R)$. 

For the exterior powers the methods of Vavilov and Petrov are unavailable. However Alexei Stepanov and the author developed a new extraction technique based on the notion of a generic element~\cite{LubStepSub}. Using Theorem~$1$ of this paper we state that the subgroup lattice $\mathcal L=L\bigl(\bw{m}\E_n(R),\GL_N(R)\bigr)$ is standard, i.\,e., the standard description of overgroups holds. Observe that for the exterior powers almost all conditions of the mentioned theorem are obvious or proved in the previous Sections. Now we must only verify the following conditions.
\begin{itemize}
\item The normilizer $N(\blank)$ is a closed subscheme in $\GL_N(\blank)$; and over a field $F$ this normalizer $N\bigl(\bw{m}\E_n(F)\bigr)$ is ``closed to be maximal'' in $\GL_N(F)$;
\item The transporter from $\bw{m}\E_n(R)$ to $N\bigl(\bw{m}\E_n(R)\bigr)$ equals $N\bigl(\bw{m}\E_n(R)\bigr)$;
\item Extract a non-trivial elementary root unipotent in a subgroup that contains the generic element of $\GL_N$ and $\bw{2}\E_n$.
\end{itemize}

Recall that the \textit{transporter} of a subgroup $E$ to a subgroup $F$ of a group $G$ is the set
$$\Tran_G(E,F)=\{g\in G\mid E^g\leq F\}.$$

It turns out that the first two items hold for an arbitrary exterior power of the elementary group. Using the invariant forms from the previous Section, we get the equality of the following four groups. Note that all normalizers and transporters here are taken in the general linear group $\GL_{N}(R)$.
\begin{theorem}
Let $R$ be a commutative ring and let $n\geq 4$. Then
$$N\bigl(\bw{m}\E_n(R)\bigr) = N\bigl(\bw{m}\SL_{n}(R)\bigr) = \Tran\bigl(\bw{m}\E_n(R), \bw{m}\SL_{n}(R)\bigr) = \bw{m}\GL_{n}(R).$$
\end{theorem}
\begin{corollary}\label{cor:TranG}
Let $R$ be a commutative ring and let $n\geq 4$. Then
$$\Tran\bigl(\bw{m}\E_n(R), \bw{m}\GL_{n}(R)\bigr) = \bw{m}\GL_{n}(R).$$
\end{corollary}
For the third item the situation is not so optimistic. The extraction is based on the decomposition of unipotens~\cite{LubReverse2}. And the latter ingredient is completely solved only for the exterior square. Thus we could extract an elementary transvection only for the exterior square too.

By summarizing the above, we get the standard description of overgroups for the exterior square under the assumptions $n\geq 6$ or $n\geq 5, 3\in R^*$. For any intermediate subgroup $H$ there exists a unique maximal ideal $A$ of the ring $R$ such that
$$\bw{2}\E_n(R)\cdot\E_N(R,A)\leq H \leq N_{\GL_N(R)}\bigl(\bw{2}\E_n(R)\cdot\E_N(R,A)\bigr).$$
Besides, the above normalizer equals the congruence subgroup of corresponding level:
$$N_{\GL_N(R)}\bigl(\bw{2}\E_n(R)\cdot\E_N(R,A)\bigr)=\rho_{A}^{-1}\left(\bw{2}\GL_{n}(R/A)\right).$$

In the general case of the $m$-th exterior powers, one can expect that under the assumption $n\geq 3m$ all intermediate subgroups $H$ are parametrized by a single ideal of the ring $R$. And conversely, if $4\leq n<3m$, then for every subgroup $H$ there is an admissible $m$-tuple of ideals $\mathbb{A}=(A_0,\dots,A_{m-1})$ in the ring $R$ such that 
$$\bw{m}\E_n(R)\cdot\E_N(R,\mathbb{A})\leq H \leq N_{\GL_N(R)}\bigl(\bw{m}\E_n(R)\cdot\E_N(R,\mathbb{A})\bigr).$$

\section{Concluding remarks}

Let us mention some further problems related to the present paper. As noted, the exterior powers of the elementary groups belong to the Aschbacher class $\mathcal{S}$. The elementary group has another irreducible representation, the symmetric power.
\begin{problem}
Describe the subgroups in $\GL_{\binom{n+m-1}{m}}(R)$ containing the symmetric power of the elementary group $\opn{S}^m\E_n(R)$ with $n\geq 3$.
\end{problem}

There is another closely related series of works in the context of linear preserver problems. For arbitrary rings the following two problems not solved. Over classically fields such as $\mathbb{C}$ or $\mathbb{R}$ some results were obtained by Vladimir Platonov, Dragomir Djokovi{\'c}, Robert Guralnick, William Waterhouse, and others, see references in the survey~\cite{VavStepSurvey}.
\begin{problem}
Obtain a description of the subgroups in $\GL_{n^2}(R)$ containing the elementary group $\E_n(R)$, $n\geq 3$ in the adjoint representation.
\end{problem}
\begin{problem}
Describe subgroups in $\GL_{n^2}(R)$ containing the elementary classical group $\EO_n(R)$, $\opn{Ep}_n(R)$, or $\opn{EU}_n(R)$ in the adjoint representation in $\GL_n(R)$.
\end{problem}


\begin{thebibliography}{XX}
\normalsize 

\bibitem{AnaVavSinI}
A.~S. Ananievsky, N.~A. Vavilov, and S.~S. Sinchuk.
\newblock {Overgroups of {$\mathrm{E}(l,R)\otimes \mathrm{E}(m,R)$} {I}.
  {L}evels and normalizers}.
\newblock {\em St. Petersbg. Math. J.}, 23(5):819--849, 2012.

\bibitem{AschbacherClasses}
M.~Aschbacher.
\newblock {On the maximal subgroups of the finite classical groups}.
\newblock {\em Invent. Math.}, 76(3):469--514, 1984.

\bibitem{BVnets}
Z.~I. Borevich and N.~A. Vavilov.
\newblock {Definition of a net subgroup}.
\newblock {\em J. Sov. Math.}, 30(1):1810--1816, 1985.

\bibitem{burness_testerman_CorrSeitz}
T.~C. Burness and D.~M. Testerman.
\newblock {Irreducible Subgroups of Simple Algebraic Groups – A Survey}.
\newblock In {\em Groups St Andrews 2017 Birmingham}, pages 230--260. Cambridge
  University Press, 2019.

\bibitem{CoopersteinStructure}
B.~N. Cooperstein.
\newblock {Nearly maximal representations for the special linear group.}
\newblock {\em Michigan Math. J.}, 27(1):3--19, 1980.

\bibitem{DieCarInvTheor}
J.~A. Dieudonn{\'{e}} and J.~B. Carrell.
\newblock {Invariant theory, old and new}.
\newblock {\em Adv. Math. (N. Y).}, 4(1):1--80, 1970.

\bibitem{KleidLiebSubStruct}
P.~B. Kleidman and M.~W. Liebeck.
\newblock {\em {The Subgroup Structure of the Finite Classical Groups}}, volume
  129.
\newblock Cambridge University Press, 1990.

\bibitem{Knus}
M.-A. Knus.
\newblock {\em {Quadratic and Hermitian Forms over Rings}}, volume 294 of {\em
  Grundlehren der mathematischen Wissenschaften}.
\newblock Springer Berlin Heidelberg, Berlin, Heidelberg, 1991.

\bibitem{LiebSeitStrucClass}
M.~W. Liebeck and G.~M. Seitz.
\newblock {On the subgroup structure of classical groups}.
\newblock {\em Invent. Math.}, 134(2):427--453, 1998.

\bibitem{LubReverse2}
R.~Lubkov.
\newblock {The reverse decomposition of unipotents for bivectors}.
\newblock {\em Commun. Algebr.}, 49(10):4546--4556, 2021.

\bibitem{LubStepSub}
R.~Lubkov and A.~Stepanov.
\newblock {Subgroups of Chevalley Groups Over Rings}.
\newblock {\em J. Math. Sci.}, 252(6):829--840, 2021.

\bibitem{SeitzMaxSub}
G.~M. Seitz.
\newblock {The maximal subgroups of classical algebraic groups}.
\newblock {\em Mem. Am. Math. Soc.}, 67(365), 1987.

\bibitem{VavSbgs}
N.~Vavilov.
\newblock {Intermediate subgroups in Chevalley groups}.
\newblock In {\em Groups Lie Type their Geom.}, pages 233--280. Cambridge
  University Press, 1995.

\bibitem{VP-EOeven}
N.~Vavilov and V.~Petrov.
\newblock {Overgroups of {$\mathrm{EO}(2l,R)$}}.
\newblock {\em J. Math. Sci.}, 116(1):2917--2925, 2003.

\bibitem{VavPere}
N.~A. Vavilov and E.~Y. Perelman.
\newblock {Polyvector representations of {$\mathrm{GL}_n$}}.
\newblock {\em J. Math. Sci.}, 145(1):4737--4750, 2007.

\bibitem{VavStepSurvey}
N.~A. Vavilov and A.~V. Stepanov.
\newblock {Overgroups of semisimple groups}.
\newblock {\em Vestn. Samar. Gos. Univ. Estestv. Ser.}, {}(3):51--95, 2008.

\bibitem{WaterhouseDet}
W.~C. Waterhouse.
\newblock {Automorphisms of {${\rm det}(X_{ij})$}: The group scheme approach}.
\newblock {\em Adv. Math. (N. Y).}, 65(2):171--203, 1987.

\end{thebibliography}

\end{document}